\documentclass[11pt]{article}
\usepackage{amssymb,epsfig,amsmath}
\usepackage[dvips]{color}
\usepackage[T1]{fontenc}
\usepackage[latin1]{inputenc}
\usepackage{comment}

\newtheorem{defi}{Definition}
\newtheorem{prop}[defi]{Proposition}
\newtheorem{theo}[defi]{Theorem}
\newtheorem{conj}[defi]{Conjecture}
\newtheorem{lemm}[defi]{Lemma}
\newtheorem{coro}[defi]{Corollary}
\newtheorem{rema}[defi]{Remark}
\newtheorem{exem}[defi]{Example}
\newtheorem{exems}[defi]{Examples}

\newcommand{\bdefi}{\begin{defi}}
\newcommand{\edefi}{\end{defi}}
\newcommand{\bprop}{\begin{prop}}
\newcommand{\eprop}{\end{prop}}
\newcommand{\btheo}{\begin{theo}}
\newcommand{\etheo}{\end{theo}}
\newcommand{\blemm}{\begin{lemm}}
\newcommand{\brema}{\begin{rema}}
\newcommand{\erema}{\end{rema}}
\newcommand{\bexer}{\begin{exem}}
\newcommand{\eexer}{\end{exem}}
\newcommand{\bexems}{\begin{exems}}
\newcommand{\eexems}{\end{exems}}
\newcommand{\bconj}{\begin{conj}}
\newcommand{\econj}{\end{conj}}
\newcommand{\elemm}{\end{lemm}}
\newcommand{\bcoro}{\begin{coro}}
\newcommand{\ecoro}{\end{coro}}
\newcommand{\dem}{\noindent{\bf Proof. }}
\newcommand{\rem}{\noindent{\bf Remark. }}

\usepackage{mathrsfs}
\renewcommand\mathcal{\mathscr}

\newcommand{\G}{{\cal G}}

\newcommand{\OOOO}{{\cal O}}

\newcommand{\HHHH}{{\cal H}}

\newcommand{\C}{{\cal C}}

\newcommand{\maths}[1]{{\mathbb #1}}  

\newcommand{\RR}{\maths{R}}
\newcommand{\NN}{\maths{N}}
\newcommand{\CC}{\maths{C}}
\newcommand{\QQ}{\maths{Q}}

\newcommand{\HH}{\maths{H}}
\newcommand{\FF}{\maths{F}}

\newcommand{\ZZ}{\maths{Z}}
\newcommand{\PP}{\maths{P}}

\newcommand{\TT}{\maths{T}}

\newcommand{\ra}{\rightarrow}
\newcommand{\bs}{\backslash}

\newcommand{\ov}[1]{{\overline #1}} 
\newcommand{\wt}[1]{{\widetilde{#1}}}
\newcommand{\wh}[1]{{\widehat{#1}}}

\newcommand{\ga}{\gamma}
\newcommand{\Ga}{\Gamma}

\newcommand{\cqfd}{\hfill$\Box$}

\newcommand{\card}{{\operatorname{Card}}}

\newcommand{\id}{\operatorname{id}}

\newcommand{\Deg}{\operatorname{deg}}
\newcommand{\Aut}{\operatorname{Aut}}

\newcommand{\hdr}{{\HH}^2_\RR}

\newcommand{\PSL}{\operatorname{PSL}}
\newcommand{\SL}{\operatorname{SL}}
\newcommand{\GL}{\operatorname{GL}}
\newcommand{\PGL}{\operatorname{PGL}}

\newcommand{\CAT}{\operatorname{CAT}}
\newcommand{\HB}{\operatorname{HB}}
\newcommand{\Par}{\operatorname{Par}}

\title{Logarithm laws for strong unstable foliations\\ 
in negative curvature\\ and non-Archimedian Diophantine approximation}
\author{Jayadev S.~Athreya \and Fr\'ed\'eric Paulin} 
\date{}

\begin{document}
\bibliographystyle{../alphanum}
\maketitle

\begin{abstract} Given for instance a finite volume negatively curved
  Riemannian manifold $M$, we give a precise relation between the
  logarithmic growth rates of the excursions into cusps neighborhoods
  of the strong unstable leaves of negatively recurrent unit vectors
  of $M$ and their linear divergence rates under the geodesic flow. As
  an application to non-Archimedian Diophantine approximation in
  positive characteristic, we relate the growth of the orbits of
  lattices under one-parameter unipotent subgroups of $\GL_2(\wh K)$
  with approximation exponents and continued fraction expansions of
  elements of the field $\wh K$ of formal Laurent series over a finite
  field. 
\end{abstract}

\section{Introduction}

The excursions of geodesic flow lines into neighborhoods of ends of
finite volume negatively curved manifolds have been studied for a long
time, and Sullivan \cite{Sullivan82} proved a seminal almost sure
logarithm law in the finite volume constant curvature case. More
probabilistic aspects have been considered too (see for instance the
works of Enriquez, Franchi, Guivarc'h, Le Jan, as
\cite{EnrFraLeJ01}). Sullivan's result has been extended by
Kleinbock-Margulis \cite{KleMar99} to finite volume locally symmetric
spaces of non compact type, by Stratmann-Velani \cite{StrVel95} to
geometrically finite constant negative curvature, by Hersonsky-Paulin
\cite{HerPau04, HerPau07} to variable negative curvature and to trees, and by
Athreya-Ghosh-Prasad \cite{AthGhoPra09,AthGhoPra12} to some
buildings.

The ergodic theory and topological dynamics of the horocyclic flow in
dimension $2$ (or of the strong unstable foliation of the geodesic
flow in higher dimension) has attracted a huge amount of studies (see
the works in negative curvature of Hedlund, Furstenberg, Dal'Bo, Dani,
Roblin, Sarig, Schapira, Smillie, as well as the recent
\cite[Chap.~9]{PauPolSha}, and in higher rank of Ratner,
Kleinbock-Margulis, Benoist-Quint and many others, see for instance
\cite{Eskin10}). But strictly analogous problems of excursions of
horocyclic flow lines or of the leaves of the strong unstable
foliation have only recently started to be studied, see for instance
the work of Athreya-Margulis \cite{AthMar09} for unipotent or
horospherical actions in some locally symmetric spaces of non compact
type, and also \cite{Athreya13,KelMoh}.

In this paper, we are interested in this problem of excursions of
(projections of) horospheres into cusps neighborhoods, with
Diophantine approximation applications, a component which is present
in all the previous works.

\medskip Let $M$ be a complete, geometrically finite, Riemannian
manifold with dimension at least $2$ and sectional curvature at most
$-1$. Let $o:T^1M\ra M$ be its unit tangent bundle and
$(\phi_t)_{t\in\RR}$ its geodesic flow. For every $v\in T^1M$,
let 
$$
W^{su}(v)= \{w\in T^1M\;:\;\lim_{t\ra+\infty}\; d(o(\phi_{-t}v),
o(\phi_{-t}w))=0\}
$$ 
be the strong unstable leaf of $v$, endowed with the naturally scaling
Hamenstadt's distance $d^{su}$ (see Section \ref{sec:geomcat}, it
coincides with the induced Riemannian distance when the sectional
curvature is constant) and, when non compact, with the filter of the
complementary subsets of its relatively compact subsets.

Our main result, Theorem \ref{theo:main}, whose simplified version is
given below, is a precise relation between the logarithmic growth
rates of the strong unstable foliation and the linear divergence rates
of the geodesic flow.

\btheo\label{theo:mainintrointro}
For every $v\in T^1 M$ which is negatively recurrent under the
geodesic flow, we have
$$  
\limsup_{w\in W^{su}(v)}\;
\frac{d(o(w),\,o(v))}{\log d^{su}(w,v)}\;=\;
  1+\limsup_{t\ra+\infty}\;\frac{d(o(\phi_{-t}v),o(v)))}{t}\;.
$$
\etheo

In the particular case when $M$ is a finite volume orientable
hyperbolic surface, we recover, by a purely geometric proof, the
logarithm law for the excursions into cusps neighborhoods of the
horocyclic flow due to \cite[Theo.~2.8]{Athreya13}.

We may also specify a set of cusps into whose neighborhoods we want to
study the excursions of the strong unstable manifolds, as follows.
Recall that a {\it cusp} of $M$ is an asymptotic class of minimizing
geodesic rays in $M$ along which the injectivity radius tends to $0$.
For every cusp $e$, let $r_e:[0,+\infty[\;\ra M$ be a representative
of $e$, and let $\beta_e:M\ra [0,+\infty[$ be the map $x\mapsto
\max\{\;0,\;\lim_{t\ra+\infty}\;t- d(x,r_e(t))\;\}$.  (One way to
normalize $\beta_e$ is to ask for $r_e$ to be contained in the closure
of, and start from the boundary of, a maximal open Margulis
neighborhood of $e$ (see for instance
\cite{BusKar81,Bowditch95,HerPau04}, it is a canonical neighborhood of
the end of $M$ to which converges $e$ if $M$ has finite volume).)
Given a (necessarily finite since $M$ is geometrically finite) set $E$
of cusps, let $\beta_E=\max_{e\in E}\beta_e$.

\btheo \label{theo:selectcuspintro}
For every $v\in T^1 M$ which is negatively recurrent under the
geodesic flow, we have
$$  
\limsup_{w\in W^{su}(v)}\;
\frac{\beta_E(o(w))}{\log d^{su}(w,v)}\;=\;
  1+\limsup_{t\ra+\infty}\;\frac{\beta_E(o(\phi_{-t}v))}{t}\;.
$$
\etheo

We refer to Corollary \ref{coro:measconseq} for almost everywhere
consequences of these theorems for the excursions of the strong
unstable leaves in cusp neighborhoods.

These theorems are valid when $M$ is replaced by the quotient of any
proper $\CAT(-1)$ metric space $X$ by any geometrically finite discrete group
of isometries of $X$, see Section \ref{sec:loglaw}.

\medskip Let us now give an application of our main result to
non-Archimedian Diophantine approximation in positive characteristic
(see for instance \cite{Lasjaunias00,Schmidt00} for nice introductions).

Let $k=\FF_q$ be a finite field with $q$ elements, where $q$ is a
positive power of a prime $p$. Let $A=k[X]$ be the ring of polynomials
in one variable $X$ over $k$ and let $K=k(X)$ be its fraction field,
endowed with the absolute value $|\cdot|=|\cdot|_\infty$ defined by
$$
\Big|\,\frac{P}{Q}\,\Big|=q^{\Deg P-\Deg Q}\;.
$$
Let $\wh K$ be the completion of $K$ for this absolute value, which is
the field $k((T^{-1}))$ of formal Laurent series $f=\sum_{i\in\ZZ}
f_iT^{-i}$ (where $f_i\in k$ is zero for $i\in\ZZ$ small enough), with
absolute value $|\cdot|=|\cdot|_\infty$ defined by
$$
|\,f\,|=q^{-\sup\{j\in\ZZ\;:\;\forall\;i<j,\;f_i=0\}}\;.
$$
For every $f\in\wh K-K$, the {\it approximation exponent} $\nu=
\nu(f)$ of $f$ is the least upper bound of the positive numbers $\nu'$ such
that there exist infinitely many elements $\frac{P}{Q}$ in $K$ (with
$P,Q$ relatively prime) such that
$$
\Big|\,f-\frac{P}{Q}\,\Big|\leq |\,Q\,|^{-\nu'}\;.
$$
Artin's {\it continued fraction expansion} of $f\in\wh K-K$ is the
sequence $(a_i=a_i(f))_{i\in\NN}$ in $A$ with $\deg a_i>0$ if $i>0$
such that
$$
f=a_0+\frac{1}{\displaystyle a_1+\frac{1}{\displaystyle 
a_2+\frac{1}{\displaystyle a_3+\frac{1}{\ddots}}}}\;.
$$
Let $\OOOO=k[[X^{-1}]]$ be the local ring of formal power series
$f=\sum_{i\in\NN} f_iT^{-i}$ (where $f_i\in k$) in $T^{-1}$ over
$k$. An {\it $\OOOO$-lattice} is a free $\OOOO$-submodule of rank $2$
in the $\wh K$-vector space $\wh K^2$. The linear action of $\GL_2(\wh
K)$ on $\wh K^2$ induces an action of $\GL_2(\wh K)$ on the set of
$\OOOO$-lattices. For every $\OOOO$-lattice $\Lambda$, let
$\Delta(\Lambda)$ 
be the unique $n\in\NN$ such that there exists $\ga\in \SL_2(A)$ and
$\lambda\in \wh K$ such that $\lambda\ga\Lambda=\OOOO\times X^{-n}\OOOO$.
For every $f\in\wh K-K$, let $(u_g=u_g(f))_{g\in\wh K}$ be the maximal
one-parameter unipotent subgroup of $\SL_2(\wh K)$ whose projective
action on the projective line $\PP_1(\wh K)=\wh K\cup\{\infty\}$ fixes
$f$.

Using the geometric approach of \cite{Paulin02} and the Bruhat-Tits
building of $(\PGL_2,\wh K)$, we have the following result, relating,
for a given irrational formal Laurent series $f$, the growth of the
orbit of any $\OOOO$-lattice under the one-parameter unipotent
subgroup of $\SL_2(\wh K)$ fixing $f$ with the approximation exponent
of $f$ and with the continued fraction expansion of $f$.

\btheo\label{theo:applicationintro} For
every $f\in\wh K-K$, we have
$$
\limsup_{|g|\ra+\infty}\;\frac{\Delta(u_g\OOOO^2)}{\log_q |g|}=2
-\frac{2}{\nu}=
1+\limsup_{n\ra+\infty}
\frac{\log|a_{n}|}{\log\big|a_{n}\prod_{i=0}^{n-1} a_i^2\,\big|}\;.
$$
\etheo

We give versions of this result for Diophantine approximation with
congruence conditions in Section \ref{sec:diophapprox}.

\medskip\noindent{\small {\it Acknowledgement: } We thank the
  Mathematisches Forchungsinstitut Oberwolfach, where this project was
  started. The first author thanks the Université Paris-Sud (Orsay)
  were this paper was continued, and G.~A.~Margulis for many useful
  and inspiring discussions. The first author also acknowledges the
  support by NSF grant DMS 1069153. The second author thanks the
  Mittag-Leffler Institute (Djursholm) where this paper was completed, and
  J.~Parkkonen for the discussion of Corollary \ref{coro:fromparpau}.}

\section{Background on $\CAT(-1)$ spaces}
\label{sec:geomcat}

We refer to \cite{BriHae99} for the definitions and basic properties
of $\CAT(-1)$ spaces, and the knowledgeable reader may skip this
section.

\medskip Let $(X,d)$ be a proper $\CAT(-1)$ geodesic metric space, and
$\ov{X}^{\rm geo}= X\cup\partial_\infty X$ its cone-topology
compactification by the asymptotic classes of its geodesic rays.

We denote by $T^1X$ the space of geodesic lines in $X$, that is, of
isometric maps $v:t\mapsto v_t$ from $\RR$ into $X$. To simplify the
notation, we will denote by $v_t$ instead of $v(t)$ the image by $v\in
T^1X$ of $t\in\RR$. When $X$ is a (complete, simply connected)
Riemannian manifold (with dimension at least $2$ and sectional
curvature at most $-1$), this notation coincides with the usual one,
upon identifying a unit tangent vector and the geodesic line it
defines. We denote by $v_\pm\in\partial_\infty X$ the points at
$\pm\infty$ of any $v\in T^1X$. The geodesic flow $(\phi_t)_{t\in\RR}$
is the action of $\RR$ on $T^1X$ by translations at the source:
$\phi_tv:s\mapsto v_{s+t}$ for all $s,t\in\RR$ and $v\in T^1X$.

\medskip
The {\it Busemann cocycle} is the continuous map
$\beta: \partial_\infty X\times X\times X\ra \RR$, defined
by
$$
(\xi,x,y)\mapsto 
\beta_\xi(x,y)=\lim_{t\ra+\infty} d(x,\xi_t)-d(y,\xi_t)\;,
$$
where $t\mapsto \xi_t$ is any geodesic ray converging to $\xi$.  For
every $\xi\in\partial_\infty X$, the {\it horospheres centered at
  $\xi$} are the level sets $f^{-1}(\lambda)$ for $\lambda\in\RR$ of
the map $f: y\mapsto \beta_\xi(y,x_0)$ from $X$ to $\RR$, and the
(closed) {\it horoballs centered at $\xi$} are its sublevel sets
$f^{-1}(]-\infty,\lambda])$ for $\lambda\in\RR$, for some (hence any)
$x_0\in X$.

If $C$ is a nonempty closed convex subset of $X$ and
$\xi\in\partial_\infty X-\partial_\infty C$, the {\it closest point to
  $\xi$ on $C$} is the unique point of $C$ which minimizes the map
$y\mapsto \beta_\xi(y,x_0)$, for some (hence any) given $x_0\in X$.

\medskip
Let $\Ga$ be a discrete group of isometries of $X$.  We denote by
$\pi:X\ra \Ga\bs X$ the canonical projection of $X$ onto its quotient
metric space $\Ga\bs X$, whose distance is again denoted by $d$.  

The limit set of $\Ga$ will be denoted by $\Lambda\Ga$, and the convex
hull of this limit set by $\C\Lambda\Ga$. Recall that $\Ga$ is
nonelementary if $\card(\Lambda\Ga)\geq 3$.  The {\it conical limit
  set} $\Lambda_c\Ga$ of $\Ga$ is the set of points
$\xi\in \partial_\infty X$ such that there exists a sequence of orbit
points of some (hence any) $x_0\in X$ under $\Ga$ converging to $\xi$
while staying at bounded distance from a geodesic ray converging to
$\xi$. The points in $\Lambda_c\Ga$ are called the {\it conical limit points}.

A point $p\in \partial_\infty X$ is a {\it bounded parabolic point}
of $\Ga$ if it is the fixed point of a parabolic element of $\Ga$ and
if its stabilizer $\Ga_p$ in $\Ga$ acts properly with compact quotient
on $\Lambda\Ga-\{p\}$.  A discrete nonelementary group of isometries
$\Ga$ of $X$ is called {\it geometrically finite} if every element of
$\Lambda\Ga$ is either a conical limit point or a bounded parabolic
 point of $\Ga$.

Let $\Par_\Ga$ be the set of fixed points of parabolic elements of
$\Ga$. If $\Ga$ is a geometrically finite group of isometries of $X$,
then (see for instance \cite{Bowditch95}) the action of $\Ga$ on
$\Par_\Ga$ has only finitely many orbits, and there exists a
$\Ga$-equivariant family $(\HB_p)_{p\in\Par_\Ga}$ of pairwise disjoint
closed horoballs, with $\HB_p$ centered at $p$, such that the quotient
$$
\Ga\bs\big(\C\Lambda\Ga-\bigcup_{p\in\Par_\Ga} \HB_p\big)
$$ 
is compact, and any geodesic ray from the boundary of $\HB_p$ to $p$
injects isometrically by the canonical projection $\pi:X\ra \Ga\bs X$.

\medskip For every $v\in T^1X$, the {\it strong unstable leaf} of $v$
is
$$
W^{\rm su}(v)=\{w\in T^1X\;:\;\lim_{t\ra+\infty} d(v_{-t},w_{-t})=
0\}\;.
$$
The set $\{w_0\;:\;w\in W^{\rm su}(v)\}$ is exactly the horosphere
centered at $v_-$ through $v_0$.  

For every $v\in T^1X$, let $d^{su}$ be {\it Hamenst\"adt's distance}
on the strong unstable leaf of $v$, defined as follows (see
\cite[Appendix]{HerPau97}, compare with \cite{Hamenstadt89}, see also
\cite[\S 2.2]{HerPau10} for a generalisation when horoballs are
replaced by arbitrary nonempty closed convex subsets): for all
$w,w'\in W^{\rm su}(v)$,
$$
d^{su}(w,w') = \lim_{t\ra+\infty}\;e^{\frac{1}{2}d(w_{-t},w'_{-t})-t}\;.
$$
This limit exists, and Hamenst\"adt's distance is a distance inducing
the original topology on $W^{\rm su}(v)$. We will denote by
$B^{su}(w,r)$ the ball of center $w$ and radius $r$ in the metric
space $(W^{su}(v),d^{su})$. For all $t\in\RR$ and $w,w'\in W^{\rm su}
(v)$, and for every isometry $\ga$ of $X$, we have $\ga
W^{su}(v)=W^{su}(\ga v)$, $\phi_t W^{su}(v)=W^{su}(\phi_t v)$,
$d^{su}(\ga w,\ga w')= d^{su}(w,w')$ and
\begin{equation}\label{eq:contracHamdist}
d^{su}(\phi_tw,\phi_tw')=e^{t}d^{su}(w,w')\;.
\end{equation}

\medskip \rem (1) When $X$ is a Riemannian manifold with constant
sectional curvature, then Hamenst\"adt's distance is the induced
Riemannian distance on the horosphere of base points of vectors of
$W^{su}(v)$ (see for instance \cite{HerPau02a}). When $X$ is a 
complex hyperbolic space $\HH^n_\CC$, then Hamenst\"adt's distance
is a multiple of Cygan's distance, see \cite[\S 3.11]{HerPau02b}).

(2) When $X$ is a metric tree, then $d^{su}(v,w)=\min\{t\in\RR\;:\;
v_{-t}=w_{-t}\}$.

(3) Here is a coarse interpretation of Hamenst\"adt's
distance. Let $\kappa>0$ be fixed. Let $\tau$ be the map defined on
the set of couples of elements of $T^1X$ in the same strong unstable
leaf, with values in $[0,+\infty[\,$, by
$$
\tau(v,w)=\min\{t\in\;\RR\;:\;d(v_{-t},w_{-t})\leq\kappa\}\;.
$$
Then is its easy to prove that there exists a constant
$c\geq 0$, depending only on $\kappa$, such that
$$
|\,\log d^{su}(v,w)-\tau(v,w)\,|\leq c\;.
$$

\medskip
Finally, we denote by $\log$ the natural logarithm, with $\log (e)=1$.

\section{Penetration in horospheres}
\label{sec:cuspincusp}

We regroup in this section the geometric lemmas concerning the
behavior of horospheres that we will need to prove our main
theorem. We refer for instance to \cite{ParPau10GT} for more
information on the penetration properties of geodesic lines in convex
subsets of $\CAT(-1)$ spaces.

Let $X$ be a $\CAT(-1)$ geodesic metric space. We will use several
times without mention the first of the following lemmas, which is well
known and follows by comparison with a geodesic triangle with an
obtuse angle in the real hyperbolic plane $\hdr$.

\blemm \label{lem:comparison} Let $x\in X$ and $y,z\in
X\cup\partial_\infty X$ be such that $x$ is the closest point to $z$
on $[x,y]$. Let $q'\in [x,z]$ and let $q$ be the intersection point of
$[y,z]$ with the sphere or horosphere centered at $z$ and passing
through $q'$. Then

\medskip\noindent
\begin{minipage}{8cm}~~~ $$
d(x,[y,z])\leq c_1=\log(1+\sqrt{2})\;,
$$
$$
d(q,q')\leq c_2=2\log\frac{1+\sqrt{3}}{2}\;.\;\;\;\Box$$
\end{minipage}
\begin{minipage}{6.9cm}
\begin{center}
\input{fig_obtuse.pstex_t}
\end{center}
\end{minipage}
\elemm

For instance by \cite[Lem.~2.9]{ParPau10GT}, for every horosphere $H$
with center $\xi,$ for every $\eta\in\partial_\infty X-\{\xi\}$, for
every $x,y\in H$ such that the geodesic rays $[x,\eta[$ and $[y,\eta[$
meet $H$ only at $x$ and $y$ respectively, we have
\begin{equation}\label{eq:parpaulong}
d(x,y)\leq 2\,c_1\;.
\end{equation}

\blemm \label{lem:c1c2} Let $x,y\in X$ be two points in a horosphere
centered at $\xi\in\partial_\infty X$, let $z$ be the closest point to
$\xi$ on $[x,y]\,$, and let $z'$ be the closest point to $y$ on
$[x,\xi[\,$. Then
$$
|\,d(x,z)-d(y,z)\,|\leq 2c_1\;\;\;{\rm and}
\;\;\;|\,d(z',x)-d(z',y)\,|\leq c_2\;.
$$  
\elemm

\medskip\noindent
\begin{minipage}{10.5cm}
  \dem This is well known, we only prove the second statement. If
  $\xi_t$ is the point at distance $t$ from $x$ on $[x,\xi[\,$, and if
  $q_t$ is the intersection with $[y,\xi_t]$ of the sphere centered at
  $\xi_t$ through $z'$, then $d(z',q_t)\leq c_2$ by Lemma
  \ref{lem:comparison} and
$$
|\,d(z',x)-d(z',y)\,|\leq |\,d(z',x)-d(q_t,y)\,|+d(z',q_t)\;,
$$
and $\lim_{t\ra+\infty}\,d(z',x)-d(q_t,y)=\beta_\xi(x,y)=0$.
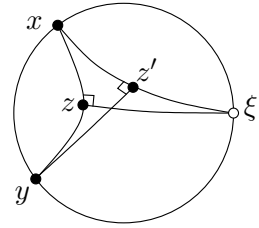
\cqfd
\end{minipage}
\begin{minipage}{4.4cm}
\begin{center}
\input{fig_piedperp.pstex_t}
\end{center}
\end{minipage}

\blemm \label{lem:controldsu}
Let $v\in T^1X$ and $w\in W^{su}(v)$. If $d^{su}(v,w)>
e^{\frac{c_2}{2}}$, then the closest point $v_{t_w}$ to $w_0$ on the
geodesic line $]v_-,v_+[$ belongs to the geodesic ray $]v_-,v_0]$, and
$$
|\,\log d^{su}(v,w)-|t_w|\,|\leq \frac{5c_1+c_2}{2}\;.
$$  
\elemm

\dem By the triangle inequality, we have
$$
d^{su}(v,w)\leq e^{\frac{1}{2}d(v_0,w_0)}\;.
$$
If $v_{t_w}$ does not belong to the geodesic ray $]v_-,v_0]$, then
$d(v_0,w_0)\leq c_2$ by Lemma \ref{lem:comparison}, which contradicts
the assumption that $d^{su}(v,w)> e^{\frac{c_2}{2}}$.  Then, using
Lemma \ref{lem:c1c2} for the last equality, writing $A=B\pm C$ instead
of $|A-B|\leq C$, we have

\medskip\noindent
\begin{minipage}{8cm}~~~ 
\begin{align*}
d(v_t,w_t)-2t&=d(v_0,w_0)\pm 4c_1\\ &
=d(v_0,v_{t_w})+d(v_{t_w},w_0)\pm 5c_1
\\ & =2|t_w|\pm(5c_1+c_2)\;.
\end{align*}
By dividing by $2$ and by taking the limit as $t$ tends to $+\infty$,
this proves the result.  \cqfd
\end{minipage}
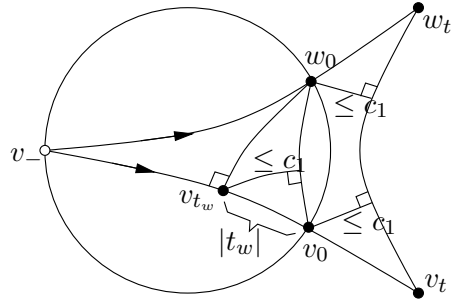
\begin{minipage}{6.9cm}
\begin{center}
\input{fig_compardsud.pstex_t}
\end{center}
\end{minipage}

\section{Horospherical logarithm laws}
\label{sec:loglaw}

Let $X$ be a proper $\CAT(-1)$ geodesic metric space, let $\Ga$ be a
geometrically finite group of isometries of $X$, and let $\pi:X\ra
\Ga\bs X$ be the canonical projection. For all $v\in T^1 X$, consider
the nondecreasing map $\Theta_v:[0,+\infty[\;\ra[0,+\infty[$ defined by
$$
\Theta_v(s)=\sup_{w\in B^{su}(v,\,s)}d(\pi(w_0),\pi(v_0))\;.
$$
A map $\psi:\;]0,+\infty[\;\ra \;]0,+\infty[$ will be called {\it
  slowly increasing} if $t\mapsto \psi(t)$ and $t\mapsto
\frac{t}{\psi(t)}$ are nondecreasing for $t$ big enough, if
$\lim_{t\ra+\infty}\psi(t)=+\infty$, and if $\lim_{t\ra+\infty}
\frac{\psi(t+c)}{\psi(t)}=1$ for all $c\in\RR$. Let $a_\psi=
\lim_{t\ra+\infty}\frac{t}{\psi(t)}\in[0,+\infty]$. For instance, for
all $a>0$ and $\alpha\in\;]0,1]$, the map $t\mapsto a\,t^\alpha$ is
slowly increasing with $a_\psi=a$ if $\alpha=1$ and $a_\psi=+\infty$
otherwise. From now on, we fix such a map $\psi$. We use the
convention that $+\infty+t=+\infty$ for all $t\in [0,+\infty]$. Note
that $\psi(t)\sim \frac{1}{a_\psi}\,t$ as $t\ra+\infty$ if $a_\psi\neq
0,+\infty$, where $f(t)\sim g(t)$ as $t\ra+\infty$ is Landau's usual
notation for $f(t)-g(t)=\operatorname{o}(g(t))$ as $t \ra \infty$.

\btheo \label{theo:main}
For all $v\in T^1 X$ such that $v_{-}$ is a conical limit point
of $\Ga$, we have
$$
\limsup_{s\ra+\infty}\;\frac{\Theta_v(s)}{\psi(\log s)}\;=\;
a_\psi+\limsup_{t\ra+\infty}\;\frac{d(\pi(v_{-t}),\pi(v_0))}{\psi(t)}\;.
$$
\etheo

\dem We start the proof by making some reductions.  Let us fix $v\in
T^1 X$ and denote by $\HHHH_v$ the horosphere with center $v_-$
through $v_0$. Let $(\HB_p)_{p\in\Par_\Ga}$ be a $\Ga$-equivariant
family of pairwise disjoint closed horoballs as in Section
\ref{sec:geomcat}, let $H_p=\partial \HB_p$ be the horosphere bounding
$\HB_p$, and let $X_{\Par}= \bigcup_{p\in\Par_\Ga}
\stackrel{\circ}{\HB}_p$ be the union of the interiors of the
horoballs $\HB_p$ for $p\in \Par_\Ga$. Let $\Delta$ be the diameter of
$\Ga\bs(\C\Lambda-X_{\Par})$.

Since $\Lambda\Ga$ has no isolated point and since $v_-\in\Lambda\Ga$,
the strong unstable leaf $W^{su}(v)$ is non compact. We endow
$W^{su}(v)$ with the filter of the complementary subsets of its
relatively compact subsets, and we will consider limits and upper
limits of functions defined on $W^{su}(v)$ along this filter.  What we
have to prove is
\begin{equation}\label{eq:maintheo}
  \limsup_{w\in W^{su}(v)}\;
\frac{d(w_0,\Ga v_0)}{\psi(\log d^{su}(w,v))}\;=\;
  a_\psi+\limsup_{t\ra+\infty}\;\frac{d(v_{-t},\Ga v_0)}{\psi(t)}\;.
\end{equation}
Since $\psi$ is slowly varying and by the triangle inequality, the
validity of this formula is unchanged if we replace $v$ by any other
given element of $W^{su}(v)$. We may hence assume that
$v_0\in\C\Lambda\Ga$ (note that $\Lambda\Ga$ is not reduced to
$\{v_-\}$). Since $v_-\in\Lambda\Ga$, the negative geodesic ray
$v_{]-\infty,\,0]}$ is therefore contained in $\C\Lambda\Ga$.  Since
$\psi$ is slowly varying, by Equation \eqref{eq:contracHamdist} and by
the triangle inequality, the validity of Equation \eqref{eq:maintheo}
is unchanged by replacing $v$ by $\phi_{-t_0}v$ for any fixed $t_0\geq
0$. Since $v_-$ is a conical limit point, we may thus assume that
$v_0\in \C\Lambda\Ga-X_{\Par}$.

Let us now introduce some more notation. 

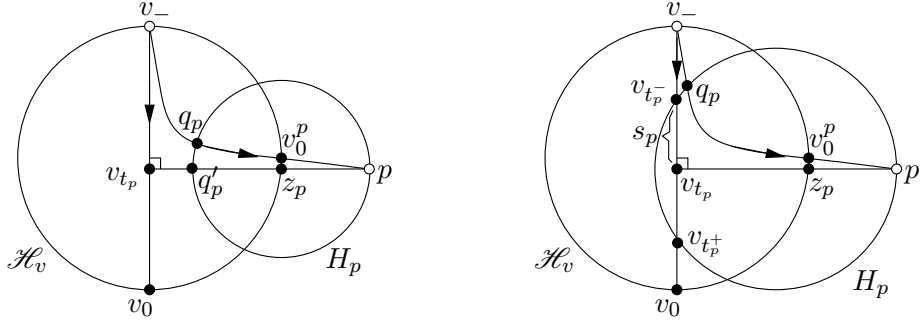
\begin{center}
\input{fig_peninp.pstex_t}
\end{center}

For every $p\in \Par_\Ga$, let $t_p\in\;]-\infty,0]$ be such that
$v_{t_p}$ is the closest point to $p$ on the geodesic ray
$v_{]-\infty,0]}$, let $z_p$ be the intersection point with $\HHHH_v$
of the geodesic ray $[v_{t_p},p[\;$, let $q_p$ be the closest point to
$v_-$ on $H_p$, and let $v^p$ be the unique element of $W^{su}(v)$
such that $v^p_+=p$. Note that the intersection with $\HHHH_v$ of the
geodesic line $v^p$ is its time $0$ point $v^p_0$. For all
$p\in\Par_\Ga$ such that $v_{t_p}\notin \HB_p$, let $q'_p$ be the
closest point to $v_{t_p}$ on $H_p$. For all $p\in\Par_\Ga$ such that
$v_{t_p}\in \HB_p$, let $t^\pm_p\in\;]-\infty,0]$ be such that
$v_{t^-_p}$ (respectively $v_{t^+_p}$) is the entering (respectively
exiting) point of the geodesic line $v$ in (respectively out) of
$\HB_p$, and let $s_p=t_p-t^-_p\geq 0$. 

By Lemma \ref{lem:controldsu}, for every $p\in\Par_\Ga$ such that
$d^{su}(v,v^p)>e^{\frac{c_2}{2}}$, we have
\begin{equation}\label{eq:controldsup}
|\,\log d^{su}(v,v^p)-|t_p|\,|\leq \frac{5c_1+c_2}{2}\;.
\end{equation}
By the initial reduction, for every $p\in\Par_\Ga$ such that
$v_{t_p}\in \HB_p$, we have $d(v_{t^-_p},\Ga x_0)\leq \Delta$,
$t^+_p\leq 0$ and $d(v_{t^+_p},\Ga x_0)\leq \Delta$. The following
estimate will also be useful.

\blemm \label{lem:controldhoro}
Let $p\in\Par_\Ga$. If $v_{t_p}\in \HB_p$, then
\begin{equation}\label{eq:controldhoro}
|\,d(v^p_0,\Ga v_0)-(|t_p|+s_p)\,|\leq 5\,c_1+c_2+\Delta\;.
\end{equation}
If $v_{t_p}\notin \HB_p$, then
\begin{equation}\label{eq:controldhorosanspene}
|\,d(v^p_0,\Ga v_0)-(|t_p|-d(v_{t_p},\HB_p))\,|\leq 
2\,c_1+2\,c_2+\Delta\;.
\end{equation}
\elemm

\dem If $v_{t_p}\in \HB_p$, since the geodesic ray $[q_p,p[$
isometrically injects in $\Ga\bs X$ and since $q_p$ belongs to
$\C\Lambda\Ga$ and is the closest point to $v^p_0$ on $H_p$, by
Equation \eqref{eq:parpaulong}, and by the second part of Lemma
\ref{lem:c1c2} for the last equality, writing $A=B\pm C$ instead of
$|A-B|\leq C$, we have
\begin{align*}
d(v_0^p,\Ga v_0)&=d(v_0^p,q_p)\pm \Delta\\ &
=d(z_p,v_{t_p^-})\pm (\Delta+4\,c_1)
\\ & =d(z_p,v_{t_p})+d(v_{t_p},v_{t_p^-})\pm (\Delta+5\,c_1)\\ 
& =|t_p|+s_p\pm (\Delta+5\,c_1+c_2)\;.
\end{align*}
The proof of the second assertion is similar: If $v_{t_p}\notin
\HB_p$, then
\begin{align*}
d(v_0^p,\Ga v_0)&=d(v_0^p,q_p)\pm \Delta\\ &=
d(z_p,q'_p)\pm(\Delta+2\,c_1+c_2)\\ &= 
d(z_p,v_{t_p})-d(v_{t_p},q'_p)\pm(\Delta+2\,c_1+c_2)\\ &= 
|t_p|-d(v_{t_p},\HB_p)\pm(\Delta+2\,c_1+2c_2)\;.\;\;\;\Box
\end{align*}

\medskip
Now that the notation is in place, let us prove Equation
\eqref{eq:maintheo} by reducing both sides to computations inside
the horoballs $\HB_p$, using the above notation.

Let us endow the set $\Par_\Ga$, and any infinite subset of it, with
the Fréchet filter of the complementary subsets of its finite
subsets. We also consider limits and upper limits of functions defined
on this set along this filter. We denote by $\stackrel{\circ}{\HB}_p$
the interior of $\HB_p$.

\blemm \label{lem:penehorodroite} If $\{p\in\Par_\Ga\;:\; v_{t_p}\in
\;\stackrel{\circ}{\HB}_p\}$ is finite, then $\limsup_{t\ra+\infty}\;
\frac{d(v_{-t},\,\Ga v_0)}{\psi(t)} =0$, and otherwise
$$
\limsup_{t\ra+\infty}\;\frac{d(v_{-t},\Ga v_0)}{\psi(t)}=
\limsup_{p\in\Par_\Ga\;:\;v_{t_p}\in\;\stackrel{\circ}{\HB}_p}\;
\frac{s_p}{\psi(|t_p|)}\;.
$$
\elemm

\dem For all $t\in\;]-\infty,0]$ such that $v_t\notin X_{\Par}=
\bigcup_{p\in\Par_\Ga}\stackrel{\circ}{\HB}_p$, we have $d(v_t,\Ga
v_0)\leq\Delta$, which in particular proves the first claim, since
$\lim_{t\ra+\infty}\psi(t)=+\infty$.  Let $p\in\Par_\Ga$ be such that
$v_{t_p}\in\;\stackrel{\circ}{\HB}_p$. For all $t\in [t^-_p,t_p]$, let
$r$ be the closest point to $v_t$ on the geodesic ray $[v_{t^-_p},p[\,$.

\medskip\noindent
\begin{minipage}{9.4cm}~~~ Since $d(v_{t_p},[v_{t^-_p},p[)\leq c_1$
  and by convexity, since the geodesic ray $[v_{t^-_p},p[$
  isometrically injects in $\Ga\bs X$ and since $v_{t^-_p}$ belongs to
  $\C\Lambda\Ga$ and is the closest point to $r$ on $H_p$, writing
  $A=B\pm C$ instead of $|A-B|\leq C$, we have
\begin{align*}
d(v_t,\Ga v_0)&=d(r,\Ga v_0)\pm c_1\\
& = d(r,v_{t^-_p})\pm (c_1+\Delta)\\ &=
d(v_t,v_{t^-_p})\pm (2c_1+\Delta)
\\ &=|t-t^-_p|\pm (2c_1+\Delta)\;.
\end{align*}
\end{minipage}
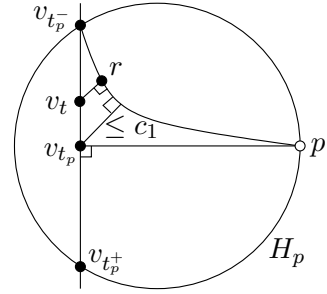
\begin{minipage}{5.5cm}
\begin{center}
\input{fig_pengeodinhoro.pstex_t}
\end{center}
\end{minipage}

\medskip\noindent Similarly, for all $t\in [t_p,t^+_p]$, we have
$|\,d(v_t,\Ga v_0)-|t-t^+_p|\;|\leq 2c_1+\Delta$.

Since $t\mapsto \psi(t)$ and $t\mapsto \frac{t}{\psi(t)}$ are
eventually nondecreasing, if $|t_p|$ is big enough, note that
$\frac{|t-t^-_p|}{\psi(|t|)}=\frac{t-t^-_p}{\psi(-t)}$ is maximal as
$t$ ranges in $]t^-_p,t_p]$ when $t=t_p$, and
$\frac{|t-t^+_p|}{\psi(|t|)} =\frac{t^+_p-t}{\psi(-t)}$ is maximal as
$t$ ranges in $[t_p,t^+_p[$ also when $t=t_p$. Since $|\,d(v_{t_p},\Ga
v_0)-s_p\;|\leq 2c_1+\Delta$, and $|\;|t_p-t^-_p|-s_p|\leq 2c_1$ by
Lemma \ref{lem:c1c2}, this proves the result.  
\cqfd

\blemm \label{lem:penehorogauche} If
$\{p\in\Par_\Ga\;:\;v_{t_p}\in\;\stackrel{\circ}{\HB}_p\}$ is finite,
then $\limsup_{w\in W^{su}(v)}\; \frac{d(w_0, \,\Ga v_0)} {\psi(\log
  d^{su}(w,\,v))}= a_\psi$ and otherwise
$$
\limsup_{w\in W^{su}(v)}\; \frac{d(w_0,\, \Ga v_0)}{\psi(\log
  d^{su}(w,\,v))}=
a_\psi+\limsup_{p\in\Par_\Ga\;:\;v_{t_p}\in\;\stackrel{\circ}{\HB}_p}\;
\frac{s_p}{\psi(|t_p|)}\;.
$$
\elemm

\dem Since $v_-\in\Lambda_c\Ga$, there exist $C>0$ and a sequence
$(\ga_n)_{n\in\NN}$ in $\Ga$ such that $\lim_{n\ra+\infty}\;
d(v_0,\ga_n v_0) =+\infty$ and $\sup_{n\in\NN}\;
d(\ga_nv_0,v_{]-\infty,\,0]})\leq C$. In particular, since the family
of pairwise disjoint horoballs $(\HB_p)_{p\in\Par_\Ga}$ is locally
finite in $X$, there exists a sequence $(p_n)_{n\in\NN}$ in $\Par_\Ga$
such that $\lim_{n\ra+\infty}\;t_{p_n}= -\infty$ and $\sup_{n\in\NN}
\;d(v_{t_{p_n}},\HB_{p_n})\leq C+\Delta$.  By Equation
\eqref{eq:controldsup} and Lemma \ref{lem:controldhoro}, since $\psi$
is eventually nondecreasing, for all $p\in\Par_\Ga$ such that
$d^{su}(v,v^p)$ is big enough, we have
$$
\frac{d(v^p_0, \Ga v_0)}{\psi(\log d^{su}(v^p,v))}\geq \min\Big\{
\frac{|t_p|-d(v_{t_p},\HB_p)-2c_1-2c_2-\Delta}
{\psi(|t_p|+\frac{5c_1+c_2}{2})},\;
\frac{|t_p|-5c_1-c_2-\Delta}{\psi(|t_p|+\frac{5c_1+c_2}{2})}
\Big\}\;.
$$
In particular, 
\begin{equation}\label{eq:minoparun}
\limsup_{w\in W^{su}(v)}\; \frac{d(w_0, \Ga v_0)}{\psi(\log
  d^{su}(w,v))}\geq \limsup_{n\ra+\infty} \frac{d(v^{p_n}_0, \Ga
  v_0)}{\psi(\log d^{su}(v^{p_n},v))}\geq a_\psi\;.
\end{equation}

Let $w\in W^{su}(v)$ with $d^{su}(w,v)> e^{\frac{c_2}{2}}$. Let
$v_{t_w}$ be the closest point to $w_+$ on the geodesic line $v$,
which belongs to $v_{]-\infty,\,0]}$ by Lemma \ref{lem:controldsu} (see
the picture below). Let $z_w$ be the intersection point with $\HHHH_v$
of the geodesic ray $[v_{t_w},p[\,$, which satisfies $d(z_w,w_0)\leq
2c_1$ by Equation \eqref{eq:parpaulong}.

If $v_{t_w}\notin X_{\Par}= \bigcup_{p\in\Par_\Ga}
\stackrel{\circ}{\HB}_p$, then respectively by the triangle
inequality, since $v_{t_w}\in \C\Lambda\Ga$, by the second claim of
Lemma \ref{lem:c1c2}, and by Lemma \ref{lem:controldsu}, we have
\begin{align*}
  d(w_0,\Ga v_0)&\leq d(z_w,\Ga v_0)+2c_1\leq
  d(z_w,v_{t_w})+2c_1+\Delta\leq |t_w|+2c_1+c_2 +\Delta\\ &\leq \log
  d^{su}(w,v)+\frac{9c_1+3c_2}{2}+\Delta\;.
\end{align*}
In particular, if $\{p\in\Par_\Ga\;:\; v_{t_p}\in \;
\stackrel{\circ}{\HB}_p\}$ is finite, since $\psi$ is slowly
incrasing, we have 
$$
\limsup_{w\in W^{su}(v)}\;
\frac{d(w_0, \Ga v_0)}{\psi(\log d^{su}(w,v))}\leq
\limsup_{w\in W^{su}(v)}\;
\frac{\log d^{su}(w,v)}{\psi(\log d^{su}(w,v))}
$$ 
is at most $a_\psi$, hence is equal to $a_\psi$ by Equation
\eqref{eq:minoparun}. This proves the first claim of Lemma
\ref{lem:penehorogauche}.

We may hence assume that $\{p\in\Par_\Ga\;:\; v_{t_p}\in \;
\stackrel{\circ}{\HB}_p\}$ is infinite. In particular,
\begin{align}
\limsup_{w\in W^{su}(v)}\; \frac{d(w_0, \Ga v_0)}{\psi(\log d^{su}(w,v))}&=
\limsup_{w\in W^{su}(v)\,:\;v_{t_w}\in\,X_{\Par}}\; 
\frac{d(w_0, \Ga v_0)}{\psi(\log  d^{su}(w,v))}\nonumber
\\ & \geq
\limsup_{p\in\Par_\Ga\;:\;v_{t_p}\in\;\stackrel{\circ}{\HB}_p}
\frac{d(v^p_0, \Ga v_0)}{\psi(\log d^{su}(v^p,v))}\;.
\label{eq:minoinfty}
\end{align}
Let us prove that the converse inequality holds. Since $\psi$ is
slowly increasing, by Equation \eqref{eq:controldsup} and the first
part of Lemma \ref{lem:controldhoro}, this will prove Lemma
\ref{lem:penehorogauche}.

Let $w\in W^{su}(v)$ and $p\in\Par_\Ga$ such that $d^{su}(w,v)>
e^{\frac{c_2}{2}}$ and $v_{t_w}\in\;\stackrel{\circ}{\HB}_p$. In
particular, $v_{t_p}\in\HB_p$. Assume for instance that
$v_{t_w}\in[v_{t_p},v_{t_p^+}[\,$.

\medskip\noindent
\begin{minipage}{9.4cm} ~~~ Let $z'_w$ be the intersection point with
  $\HHHH_v$ of the geodesic ray $[v_{t^+_p},p[\,$, which satisfies
  $d(z'_w,w_0)\leq 2c_1$ by Equation \eqref{eq:parpaulong}. Then, again
  using the second claim of Lemma \ref{lem:c1c2} for the final inequality,
\begin{align*}
d(w_0,\Ga v_0)&\leq d(z'_w,\Ga v_0)+2c_1\leq
 d(z'_w,v_{t^+_p})+2c_1+\Delta\\ &\leq
 d(z_w,v_{t^+_p})+4c_1+\Delta\\ &\leq 
d(z_w,v_{t_w})+d(v_{t_w},v_{t^+_p})+5c_1+\Delta
\\ &\leq |t_w|+|t_w-t^+_p|+5c_1+c_2+\Delta\;.
\end{align*}
\end{minipage}
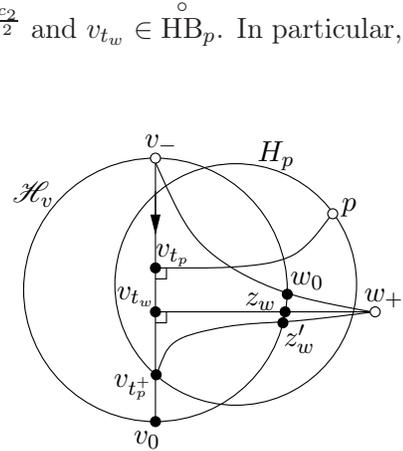
\begin{minipage}{5.5cm}
\begin{center}
\input{fig_controlw.pstex_t}
\end{center}
\end{minipage}

\medskip\noindent The map from $[t_p,t_p^+[$ to $[0,+\infty[$ defined
by $t_w\mapsto \frac{|t_w|+|t_w-t^+_p|}{\psi(|t_w|)}=
\frac{t^+_p-2t_w} {\psi(-t_w)}$ is nonincreasing if $t^+_p$ is small
enough, with maximum reached at $t_w=t_p$. This maximum is at least
$\frac{|t_p|+s_p- c_2} {\psi(|t_p|)}$ by the first assertion of Lemma
\ref{lem:c1c2}.  Since $\lim_{t\ra+\infty}\psi(t)=+\infty$, this
proves the converse part of Equation \eqref{eq:minoinfty}, thus proves
Lemma \ref{lem:penehorogauche}.  
\cqfd

\medskip
Now Equation \eqref{eq:maintheo}, hence  Theorem
\ref{theo:main}, follows immediately from Lemma
\ref{lem:penehorodroite} and Lemma \ref{lem:penehorogauche}.  \cqfd

\medskip For every $p\in\Par_\Ga$, let $\wt\beta_p:X\ra[0,+\infty[$ be
the (well-defined and $1$-Lipschitz) map
$$
\wt\beta_p:x\mapsto\max\{0,\lim_{t\ra+\infty}t-d(x,r_p(t))\}\;,
$$
where $r_p:[0,+\infty[\;\ra X$ is any geodesic ray from a point of $H_p$
to $p$.
For every $\Ga$-invariant subset $E$ of $\Par_\Ga$, let
$$
\wt\beta_E=\max_{p\in E}\;\wt \beta_p\;,
$$
which is a $\Ga$-invariant $1$-Lipschitz map from $X$ to
$[0,+\infty[\,$. The proof of the following result is the same as the
one of Theorem \ref{theo:main}, up to replacing the full family
$(\HB_p)_{p\in\Par_\Ga}$ by the subfamily $(\HB_p)_{p\in E}$.

\btheo \label{theo:selectcusp}
For every $\Ga$-invariant subset $E$ of $\Par_\Ga$, for every
$v\in T^1 X$ such that $v_{-}\in\Lambda_c\Ga$, we have
$$  
\limsup_{w\in W^{su}(v)}\;
\frac{\wt\beta_E(w_0)}{\psi(\log d^{su}(w,v))}\;=\;
  a_\psi+\limsup_{t\ra+\infty}\;\frac{\wt\beta_E(v_{-t})}{\psi(t)}\;.\; \;\;\Box
$$
\etheo

Theorem \ref{theo:selectcuspintro} in the introduction is a corollary
of Theorem \ref{theo:selectcusp}, since replacing $\wt \beta_p$ by
$\wt \beta_p+c_p$ for any constant $c_p\in\RR$ depending only on the
orbit of $p$ under $\Ga$ does not change its validity.

\medskip Note that under the hypothesis of an almost sure logarithm law
for an invariant measure of the geodesic flow, the right hand sides
of the equations in Theorem \ref{theo:main} and Theorem
\ref{theo:selectcusp} are an easily computed constant for almost every
$v$, hence so are the left hand sides. In particular, the following
result follows from \cite[Coro.~6.1]{HerPau04} and
\cite[Coro.~1.2]{HerPau07}. The new assumptions are satisfied in
particular if the Riemannian metric of $X$ is locally symmetric in at
least one horoball centered at each parabolic point. We refer for
instance to \cite{Roblin03} for the definitions and properties of the
critical exponents and of the Bowen-Margulis measure.

\bcoro \label{coro:measconseq} Let $(X,\Ga)$ be as above, with
$\Par_\Ga$ nonempty. Let $\delta$ be the critical exponent of
$\Ga$. Assume furthermore that either $X$ is a locally finite tree, or
$X$ is a Riemannian manifold with pinched sectional curvature such
that, for every $p\in\Par_\Ga$, if $\delta_p$ is the critical exponent
of the stabilizer $\Ga_p$ of $p$ in $\Ga$, then $\delta_p<\delta$ and
there exists $c>0$ such that $\frac{1}{c}\,e^{\delta_p n}\leq
\card\{\alpha\in\Ga_p\;:\; d(x_0,\alpha x_0)\leq n\}\leq c\;e^{\delta_
  p n}$ for all $n\in\NN$.

Then for almost every $v\in T^1X$ for the Bowen-Margulis measure of
$\Ga$ on $T^1X$, we have 
$$
\limsup_{s\ra+\infty}\;
\frac{\Theta_v(s)}{\psi(\log s)}=a_\psi\;.
$$
\ecoro

\medskip \dem Let $\delta_0=\max_{p\in\Par_\Ga} \delta_p<\delta$.

First note that, by \cite[\S 2]{Roblin03}, the set of elements $v\in
T^1X$ such that $v_+\in\Lambda_c\Ga$ has full measure for the
Bowen-Margulis measure. Also note that since $\psi$ is slowly varying,
if $a_\psi\neq+\infty$, then $\lim_{t\ra+\infty}\frac{\log t}{\psi(t)}
=0$. Hence if $a_\psi\neq+\infty$, by Theorem \ref{theo:main}, we
have, for almost every $v\in T^1X$ for the Bowen-Margulis measure, if
$\limsup_{t\ra+\infty}\;\frac{d(v_{-t},\,\Ga v_0)}{\log t}$ is finite,
then
\begin{equation}\label{eq:reducAthPautoHerPau}
\limsup_{s\ra+\infty}\;\frac{\Theta_v(s)}{\psi(\log s)}\;=\;a_\psi+
\limsup_{t\ra+\infty}\;\frac{d(v_{-t},\Ga v_0)}{\log t}
\lim_{t\ra+\infty}\frac{\log t}{\psi(t)}=a_\psi\;,
\end{equation}
a formula which is also true if $a_\psi=+\infty$.

In the locally finite tree case, \cite[Coro.~1.2]{HerPau07} applies
directly, since it proves that $\limsup_{t\ra+\infty}\;
\frac{d(v_{-t},\Ga v_0)}{\log t}=\frac{1}{\delta}$ for almost every
$v\in T^1X$ for the Bowen-Margulis measure.

Assume hence that $X$ is a Riemannian manifold as in the
statement. The only assumption of \cite[Coro.~6.1]{HerPau04} that is
not an assumption of Corollary \ref{coro:measconseq} is that there
exists $c'>0$ such that $\frac{1}{c'}\,e^{\delta n}\leq
\card\{\ga\in\Ga\;:\; d(x_0,\ga x_0)\leq n\}\leq c'\;e^{\delta n}$ for
all $n\in\NN$.

By \cite{DalOtaPei00}, the assumptions that $\Ga$ is geometrically
finite and that $\delta_p<\delta$ for all $p\in\Par_\Ga$ imply the
finiteness of the Bowen-Margulis measure on $\Ga\bs T^1X$.  By
\cite{DalBo99} since $\Par_\Ga\neq \emptyset$, the set of the
translation lengths of the hyperbolic elements of $\Ga$ is not
contained in a discrete subgroup of $\RR$.  By \cite{Roblin03}, the
extra assumption above is satisfied (and there is even an asymptotic
equivalent $\card\{\ga\in\Ga\;:\; d(x_0,\ga x_0)\leq n\}\sim
c'\;e^{\delta n}$ as $n\ra+\infty$).  Since the conclusion of
\cite[Coro.~6.1]{HerPau04} is that $\limsup_{t\ra+\infty}
\;\frac{d(v_{-t},\,\Ga v_0)}{\log t}= \frac{1}{2(\delta-\delta_0)}$ for
almost every $v\in T^1X$ for the Bowen-Margulis measure, Corollary
\ref{coro:measconseq} follows from Equation
\eqref{eq:reducAthPautoHerPau}.  \cqfd

\medskip \rem Since the results \cite[Coro.~6.1]{HerPau04} and
\cite[Coro.~1.2]{HerPau07} are valid cusp by cusp, a statement
analogous to Corollary \ref{coro:measconseq} for a prescribed set of
cusps is also valid.

\medskip
As another application of our main theorem, here is another consequence,
for the behavior of strong unstable leaves, of properties of the
geodesic flow.

\bcoro \label{coro:fromparpau}
If $\Ga$ is convex-cocompact, then for every $v\in T^1X$ such
that $v_-\in\Lambda\Ga$, we have $\limsup_{s\ra+\infty}\;
\frac{\Theta_v(s)}{\log s}=1$.

If $\Ga$ is not convex-cocompact and if $X$ is a Riemannian manifold
of dimension at least $3$, then for every $\alpha\in[1,2]$, there
exists $v\in T^1X$ such that $\limsup_{s\ra+\infty}\;
\frac{\Theta_v(s)}{\log s}=\alpha$.  
\ecoro

\dem The first claim is immediate from Theorem
\ref{theo:mainintrointro}.  The second one follows from the techniques
of \cite[\S 5.4]{ParPau10GT}.  \cqfd

\medskip
Given $\alpha\in \;]1,2]$, it would be interesting to study the
Hausdorff dimension of the set of elements $v\in T^1X$ such that 
$\limsup_{s\ra+\infty}\;\frac{\Theta_v(s)}{\log s}=\alpha$.

\section{An application to non-Archimedian Diophantine 
approximation}
\label{sec:diophapprox}

For all $n\geq 2$, let $\HH^n_\RR$ be the upper halfspace model of the
real hyperbolic space of dimension $n$. Applications to Archimedian
Diophantine approximation may be obtained, as in the case of
$X=\HH^2_\RR$ and $\Ga$ a congruence subgroup of $\PSL_2(\ZZ)$ (see
for instance \cite{AthMar09}), by taking for instance $X=\HH^3_\RR$
and $\Ga=\PSL(\OOOO)$ where $\OOOO$ is an order in the ring of
integers of an imaginary quadratic number field, or $X=\HH^5_\RR$ and
$\Ga=\PSL(\OOOO)$ where $\OOOO$ is an order in a definite quaternion
algebra over $\QQ$ (see for instance \cite{parPau12ANT}). But in this
paper, we concentrate on the applications to non-Archimedian
Diophantine approximation.

\medskip We start this section by restating a version of Theorem
\ref{theo:selectcusp} in the particular case of trees, which will be
more directly applicable for our arithmetic applications.

Let $T$ be a localy finite tree (endowed with the maximal distance
making each edge isometric to $[0,1]$, which is $\CAT(-1)$). Let $VT$
be its set of vertices and $\Aut(T)$ its locally compact automorphism
groups (which is contained in its isometry group). Let $\Ga$ be a
geometrically finite subgroup of $\Aut(T)$. Up to taking the first
barycentric subdivision of $T$ and rescaling, we assume that $\Ga$
acts without inversion (that is, no element of $\Ga$ maps an edge of
$T$ to its opposite edge), so that $\Ga\bs T$ has a unique structure
of graph such that the canonical projection $T\ra \Ga\bs T$ is a
morphism of graphs. By the structure theorem of \cite{Paulin04b}
(improving on the algebraic cases of Serre \cite{Serre83} and Lubotzky
\cite{Lubotzky91}), with $E$ the finite set $\Ga\bs\Par_\Ga$, there
exist a finite subgraph $\G$ of $\Ga\bs\C\Lambda\Ga$ and for every
$e\in E$, a geodesic ray $\rho_e:[0,+\infty[\;\ra \Ga\bs\C\Lambda\Ga$
with origin a vertex, which lifts to a geodesic ray in $T$ converging
to any representative of $e$ in $\Par_\Ga$, such that $\Ga\bs
\C\Lambda\Ga$ is the disjoint union of $\G$ and the open rays
$\rho_e(]0,+\infty[)$ for $e\in E$.

For every $e\in E$, define a map $\Delta_e:VT\ra [0,+\infty[$ by
$\Delta_e(x)=n$ if $\Ga x= \rho_e(n)$ (such an $n$ is unique if it
exists), and $\Delta_e(x)=0$ otherwise. Note that if $\wt e$ is an
element of $\Par_\Ga$ whose image in $\Ga\bs\Par_\Ga$ is $e$, since
$\rho_e$ lifts to a geodesic ray converging to $\wt e$, there exists a
constant $c\in\RR$ such that $\wt \beta_{\{\wt e\}}(t)=\Delta_e(t)+c$
for $t$ big enough, with the notation before Theorem
\ref{theo:selectcusp}.  Also note that two geodesic lines in $T$,
starting from the same point at infinity, coincide up to translation
on a neighborhood of $-\infty$.  For every non-isolated point
$\xi_*\in\partial_\infty T$, we endow $\partial_\infty T -\{\xi_*\}$
with the filter of the complementary subsets of its relatively compact
subsets. Therefore, the following result follows immediately from the
definition of Hamenstädt's distance and Theorem \ref{theo:selectcusp}.

\bcoro \label{coro:applicationtreecase} Let $\xi_*\in\Lambda_c\Ga$ and
$\eta_*\in\partial_\infty T-\{\xi_*\}$. For every
$\eta\in\partial_\infty T-\{\xi_*\}$, let $t\mapsto \eta(t)$ be the
geodesic line from $\xi_*$ to $\eta$ such that $\eta_*(0)\in VT$ and
$\eta(t)=\eta_*(t)$ for $t$ small enough. Let
$\delta_*(\eta,\eta_*)=\inf\{t\in\NN\;:\; \forall\; s\geq t,\;
\eta(-t)=\eta_*(-t)\}$. Then for all $e\in E$ and
$\psi:\;]0,+\infty[\;\ra\;]0,+\infty[$ slowly increasing, we have
$$
\limsup_{\eta\in \partial_\infty T -\{\xi_*\}}\;
\frac{\Delta_e(\eta(0))}{\psi(\delta_*(\eta,\eta_*))} =
a_\psi+\limsup_{t\ra+\infty}\;
\frac{\Delta_e(\eta_*(-t))}{\psi(t)}\;. \;\;\;\Box
$$
\ecoro

\medskip Let us now give our applications to non-Archimedian
Diophantine approximation. We follow the notation of \cite{Paulin02},
in particular as recalled in the introduction for $k=\FF_q$, $A=k[X]$,
$K=k(X)$, $\wh K=k((X^{-1}))$, $\OOOO=k[[X^{-1}]]$,
$|\cdot|=|\cdot|_\infty$ and, for every $f\in\wh K-K$, its continued
fraction expansion $(a_n= a_n(f))_{n\in\NN}$ and its approximation
exponent $\nu=\nu(f)$. 

\medskip For every $f=\sum_{i\in\ZZ} f_iT^{-i}\in\wh K$, the {\it
  integral part} $[f]$ of $f$ is $\sum_{i\leq 0} f_iT^{-i}\in A$ and
its {\it fractional part} $\{f\}$ is $\sum_{i>0} f_iT^{-i}\in
X^{-1}\OOOO$.  {\it Artin's map} $\Psi: X^{-1}\OOOO-\{0\}\ra
X^{-1}\OOOO$ is defined by $f\mapsto \{1/f\}$. Given $f\in\wh K-K$, we
have $a_0=[f]$ and if $n\geq 1$, then
$a_n=[\frac{1}{\psi^{n-1}(f-a_0)}]$. Consider the sequences
$(P_n)_{n\in\NN\cup\{-1\}}$ and $(Q_n)_{n\in\NN\cup\{-1\}}$ in $A$
inductively defined by
$$
P_{-1}=1, Q_{-1}=0, P_0=a_0, Q_0=1
$$
and for every $n\in\NN$
$$
P_{n+1}=a_{n+1}P_n+P_{n-1}\;\;\;{\rm and}\;\;\;
Q_{n+1}=a_{n+1}Q_n+Q_{n-1}\;.
$$ 
Then $P_n$ and $Q_n$ are relatively prime, and
$$
\frac{P_n}{Q_n}=a_0+\frac{1}{\displaystyle a_1+\frac{1}{\displaystyle 
a_2+\frac{1}{\displaystyle \begin{array}{c} \ddots 
\;\;\;\;\;\;\;\;\;\;\;\; \\
{\displaystyle a_{n-1}+\frac{1}{a_n}}\end{array}}}}
$$
is called the {\it $n$-th convergent} of $f$. The sequence
$(\frac{P_n}{Q_n})_{n\in\NN}$ converges to $f$ (for the above, see for
instance \cite{Lasjaunias00,Schmidt00}, as well as \cite{Paulin02} for
a geometric explanation).

The action of $\GL_2(\wh K)$ on the set of $\OOOO$-lattices induces an
action of $\wh K^*$ by homotheties on this set, and we will denote by
$[\Lambda]$ the homothety class of an $\OOOO$-lattice $\Lambda$.

\medskip
\noindent {\bf Remarks. }
(1) Note that $\OOOO$-lattices $\Lambda$ in $\wh K^2$ behave, from the
topological viewpoint, very differently than $\ZZ$-lattices in
$\RR^2$: they are compact open additive subgroups of $\wh K^2$, and
hence $\wh K^2/\Lambda$ is infinite and discrete (thus non
compact). Furthermore, for any norm $\|\cdot\|$ on $\wh K^2$, we have
$\inf_{x\in\Lambda-\{0\}}\|x\| =0$.

(2) The set $V\TT_{\wt K}$ (see below for an explanation of this
notation) of homothety classes of $\OOOO$-lattices in $\wh K^2$ can be
endowed with the quotient of Chabauty's topology on closed subgroups
of the (additive) locally compact group $\wh K^2$, or, equivalently,
with the topology of an homogeneous space under the transitive
(linear) action of $\PGL_2(\wh K)$. Note that this topology is
discrete since $\PGL_2(\OOOO)$ is open in $\PGL_2(\wh K)$, again a
major difference from the case of $\ZZ$-lattices in $\RR^2$. The map
$\Delta$ defined in the introduction induces (by passing to the
quotients) a proper map from $\PSL_2(A)\bs V\TT_{\wt K}$ to
$\NN$. This map is an ultrametric analog of the inverse of the systole
map on $\ZZ$-lattices with covolume $1$ in $\RR^2$, whose properness
is called Mahler's criterion.

\medskip
Let us fix a nonzero element $Q_*$ of $A$. Consider Hecke's
nonprincipal congruence subgroup
$$
\Ga^0_{Q_*}=\Big\{\Big(\begin{array}{cc} a & b \\ c & d\end{array}
\Big)\in\SL_2(A)\;:\; c\equiv 0\mod Q_*\Big\}\;,
$$
which has finite index in $\SL_2(A)$. For every $\OOOO$-lattice
$\Lambda$, define $\Delta_{Q_*}(\Lambda)=n$ if there exists $\ga\in
\Ga^0_{Q_*}$ such that $[\ga \Lambda]= [\OOOO\times X^{-n}\OOOO]$ and
$\Delta_{Q_*}(\Lambda)=0$ otherwise. 

For every $f\in\wh K-K$, recall that $(u_g=u_g(f))_{g\in\wh K}$ is the
maximal one-parameter unipotent subgroup of $\SL_2(\wh K)$ whose
projective action on $\PP_1(\wh K)=\wh K\cup\{\infty\}$ fixes $f$. For
every $f\in\wh K-K$, the {\it approximation exponent $\nu_{Q_*}=
  \nu_{Q_*}(f)$ of $f$ relative to $Q_*$} is the least upper bound of
the positive numbers $\nu'$ such that there exist infinitely many
elements $\frac{P}{Q}$ in $K$ with $P$ and $Q$ relatively prime and
$Q\equiv 0\mod Q_*$ such that
$$
\Big|\,f-\frac{P}{Q}\,\Big|\leq |\,Q\,|^{-\nu'}\;.
$$

\btheo\label{theo:application} For every $f\in\wh K-K$, we have
$$
\limsup_{|g|\ra+\infty}\;\frac{\Delta_{Q_*}(u_g\OOOO^2)}{\log_q |g|}=2
-\frac{2}{\nu_{Q_*}}= 1+\limsup_{n\ra+\infty\;:\;Q_n\equiv 0\!\!\!\mod Q_*}
\frac{\log|a_{n+1}|}{\log\big|a_{n+1}\prod_{i=1}^n a_i^2\,\big|}\;.
$$
\etheo

Theorem \ref{theo:applicationintro} in the introduction follows by
taking $Q_*=1$. 

\medskip \dem We will apply Corollary \ref{coro:applicationtreecase}
with $T$ the Bruhat-Tits tree of $(\PGL_2,\wh K)$, whose definition
and useful properties we start by recalling (following \cite{Serre83}).

The {\it Bruhat-Tits tree} $\TT_{\wh K}$ of $(\PGL_2,\wh K)$ is the
graph whose vertices are the homothety classes of $\OOOO$-lattices in
$\wh K^2$, two vertices $x$ and $x'$ being joined by an edge if and
only if there exist representatives $\Lambda,\Lambda'$ of $x,x'$
respectively such that $\Lambda'\subset \Lambda$ and
$\Lambda/\Lambda'$ is isomorphic to $\OOOO/X^{-1}\OOOO$.

We identify as usual the projective line  $\PP_1(\wh K)$ with $\wh
K\cup\{\infty\}$ by the map $\wh K^*(x,y)\mapsto \frac{x}{y}$.  We
denote by $(g,x)\mapsto g\cdot x$ the projective action of $g\in
\GL_2(\wh K)$ on $x\in\PP_1(\wh K)=\wh K\cup\{\infty\}$.

The action of $\GL_2(\wh K)$ on the set of $\OOOO$-lattices induces an
isometric action of $\GL_2(\wh K)$ on $\TT_{\wh K}$. Note that
$\SL_2(\wh K)$ acts with two orbits on the set of vertices of
$\TT_{\wh K}$. There exists one and only one homeomorphism between
$\partial_\infty \TT_{\wh K}$ and $\PP_1(\wh K)$ such that the
(continuous) extension to $\partial_\infty \TT_{\wh K}$ of the
isometric action of $\GL_2(\wh K)$ on $\TT_{\wh K}$ corresponds to the
projective action of $\GL_2(\wh K)$ on $\PP_1(\wh K)$. From now on, we
identify $\partial_\infty \TT_{\wh K}$ and $\PP_1(\wh K)$ by this
homeomorphism.

We denote by $\HB_\infty$ the horoball in $\TT_{\wh K}$ with center
$\infty$ whose boundary contains the vertex $[\OOOO^2]$.  Note that
$\SL_2(A)$ is a geometrically finite group of isometries of $\TT_{\wh
  K}$, with only one orbit of parabolic points, and that
$(\ga\HB_\infty)_{\ga\in\SL_2(A)/\SL_2(A)_\infty}$ is the associated
$\SL_2(A)$-equivariant family of maximal horoballs with pairwise
disjoint interiors (see \cite[\S 6.2]{Paulin02}.  The geodesic ray
from $[\OOOO^2]$ to $\infty\in\partial_\infty\TT_{\wh K}$, whose
sequence of consecutive vertices is $([\OOOO\times
X^{-n}\OOOO])_{n\in\NN}$, injects onto the quotient
$\SL_2(A)\bs\TT_{\wh K}$.

\medskip Let us fix $f\in\wh K-K$. We are now going to apply Corollary
\ref{coro:applicationtreecase} with $\psi=\id$ (so that $a_\psi=1$),
$T=\TT_{\wh K}$, $\Ga$ the image of $\Ga^0_{Q_*}$ in $\Aut(T)$ (which is
also a geometrically finite subgroup with $\Par_\Ga= \Par_{\PSL_2(A)}=
\PP_1(K)=K\cup\{\infty\}$), $\xi_*=f$ (which is a conical limit point,
since $f$ is irrational and the limit set of $\PSL_2(A)$, hence of
$\Ga$, is the whole boundary at infinity), $\eta_*=0$, and
$e=\Ga \infty$.

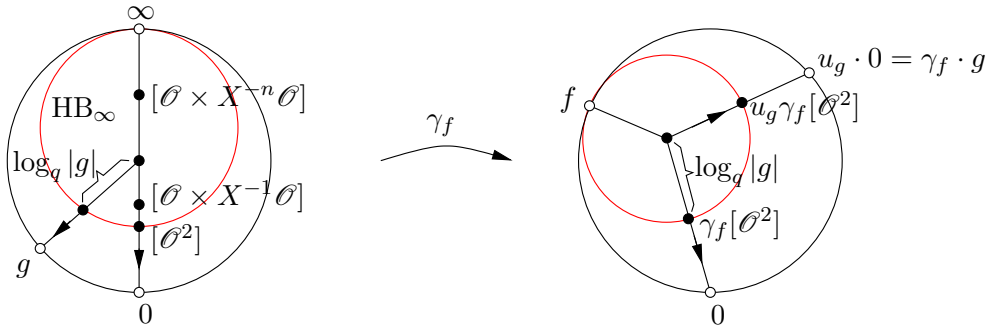
\begin{center}
\input{fig_bruhattits.pstex_t}
\end{center}

The sequence $([\OOOO\times X^{-n}\OOOO])_{n\in\NN}$ of consecutive
vertices of the geodesic ray in $\TT_{\wh K}$ from $[\OOOO^2]$ to
$\infty$ isometrically injects in $\SL_2(A)\bs\TT_{\wh K}$, hence in
$\Ga^0_{Q_*}\bs\TT_{\wh K}$. Therefore, there exists a constant
$c\in\RR$ such that for every $\OOOO$-lattice $\Lambda$, we have
$$
|\,\Delta_e([\Lambda])-\Delta_{Q_*}(\Lambda)\,|\leq c\;.
$$

The element $\ga_f=\Big(\!\begin{array}{cc}
  1 & 0\\ 1/f & 1\end{array}\!  \Big)\in \SL_2(\wh K)$ maps projectively
$\infty$ to $f$ (hence sends the horospheres centered at $\infty$ to
the horospheres centered at $f$) and fixes $0$.  The maximal
one-parameter unipotent subgroup of $\SL_2(\wh K)$ fixing $\infty$ is
$g\mapsto \Big(\!\begin{array}{cc} 1 & g \\ 0 & 1\end{array}\!
\Big)$. Hence the maximal one-parameter unipotent subgroup
$(u_g)_{g\in\wh K}$ of $\SL_2(\wh K)$ fixing $f$ is $g\mapsto u_g=
\ga_f\Big(\!\!\begin{array}{cc} 1 & g \\ 0 & 1\end{array}\!
\Big)\ga_f^{-1}$. We fix the parametrization of the geodesic lines
starting from $f$ so that they cross at time $t=0$ through
$\ga_f \partial\HB_\infty$. Note that there exists a constant
$c'\in\RR$ such that for all $g\in \wh K$, we have
$$
\big |\,\Delta_{Q_*}(u_g\ga_f\OOOO^2)-
\Delta_{Q_*}(u_g\OOOO^2)\,\big |\leq c'\;.
$$

By \cite[Coro.~5.2]{Paulin02}, for all $g\in\wh K$, Hamenstädt's
distance between the geodesic lines starting from $\infty$ and ending
at $0$ and at $g$, passing through $\partial\HB_\infty$ at time $0$,
is $|g|^\frac{1}{\log q}$. Since $\ga_f$ is an isometry and $\ga_f\dot
g= u_g\cdot 0$, with $\delta_*(\cdot,\cdot)$ defined in Corollary
\ref{coro:applicationtreecase}, we hence have
$$
\delta_*(u_g\cdot 0,0)=\log_q|g|\;.
$$ 
The map from $\wh K$ to $\partial_\infty\TT_{\wh K}-\{f\}$ sending $g$
to $\eta=u_g\cdot 0$ is a homeomorphism, and
$\eta(0)=u_g\ga_f[\OOOO^2]$. Hence
\begin{equation}\label{eq:premiermemb}
\limsup_{\eta\in \partial_\infty T -\{\xi_*\}}\;
\frac{\Delta_e(\eta(0))}{\delta_*(\eta,\eta_*)} =
\limsup_{|g|\ra+\infty}\;\frac{\Delta_{Q_*}(u_g\OOOO^2)}{\log_q |g|}\;.
\end{equation}

Let us denote by $v_\infty:\wh K\ra\ZZ\cup\{+\infty\}$ the valuation
associated to the absolute value $|\cdot|$, so that for all $f\in\wh
K$, we have
$$
|\,f\,|=q^{-v_\infty(f)}\;.
$$

By \cite[\S 6.3]{Paulin02}, the geodesic line starting from $\infty$
and ending at $f$ (passing at time $t=0$ through $\partial\HB_\infty$)
enters successively the interiors of the horoballs of the
$\SL_2(A)$-equivariant family of maximal horoballs with pairwise
disjoint interiors which are centered at the convergents of
$f$. Furthermore, if $x_n$ is its entering vertex in
$\HB_{\frac{P_n}{Q_n}}$ and $y_n$ its point the closest to
$\frac{P_n}{Q_n}$, then
$$
d(x_{n},y_{n})=\frac{1}{2}\,d(x_{n},x_{n+1})=-v_\infty(a_{n+1})\;.
$$
In particular, 
$$
d(x_0,y_n)=\sum_{i=1}^{n}-2\,v_\infty(a_{i})
-v_\infty(a_{n+1}) = \log_q\big|a_{n+1}\prod_{i=1}^n a_i^2\,\big|\;.
$$ 
By the definition of $\Delta_{Q_*}$, we are only interested in the
penetration of the geodesic line $]\infty,f[$ in the horospheres
$\HB_{\frac{P}{Q}}$ with $P,Q\in A$ relatively prime and $Q=0\mod
Q_*$. Since the geodesic ray from $\ga_f[\OOOO^2]$ to $f$ coincides,
for times big enough and up to translation, with the geodesic line
$]\infty,f[\,$, we have
\begin{equation}\label{eq:secondmemb}
\limsup_{t\ra+\infty}\;\frac{\Delta_e(\eta_*(-t))}{t}
=\limsup_{n\ra+\infty\;:\;Q_n=0\!\!\!\mod Q_*}
\frac{\log_q|a_{n+1}|}{\log_q\big|a_{n+1}\prod_{i=1}^n a_i^2\,\big|}\;.
\end{equation}

By \cite[\S 6.3]{Paulin02}, with the above notation, we also have
$$
d(x_0,y_n)=v_\infty\big(f-\frac{P_n}{Q_n}\big)=
-\log_q\big|\,f-\frac{P_n}{Q_n}\,\big|\;,
$$
and
$$
d(x_0,x_n)=-2\,v_\infty(Q_n)=2\log_q|\,Q_n\,|\;.
$$ 

Since $|\,f-\frac{P}{Q}\,|$ is at least $1$ if $\frac{P}{Q}$ is not a
convergent of $f$, we hence have, by the definition of the
approximation exponent $\nu_{Q_*}$,
\begin{align}
\frac{2}{\nu_{Q_*}}&=\liminf_{|Q|\ra+\infty\;:\;Q\equiv 0\!\!\!\mod Q_*}
\frac{2\log |\,Q\,|}{-\log|\,f-\frac{P}{Q}\,|}= 
\liminf_{n\ra+\infty\;:\;Q_n\equiv 0\!\!\!\mod Q_*}
\;\frac{2\log_q |\,Q_n\,|}{-\log_q|\,f-\frac{P_n}{Q_n}\,|}
\nonumber \\ &= \liminf_{n\ra+\infty\;:\;Q_n\equiv 0\!\!\!\mod Q_*}\;
\frac{d(x_0,x_n)}{d(x_0,y_n)}\nonumber \\ &=
\liminf_{n\ra+\infty\;:\;Q_n\equiv 0\!\!\!\mod Q_*}\;
\frac{\sum_{i=1}^{n}-2v_\infty(a_{i})}
{\sum_{i=1}^{n}-2v_\infty(a_{i})-v_\infty(a_{n+1})}\nonumber \\ &=
1-\limsup_{n\ra+\infty\;:\;Q_n\equiv 0\!\!\!\mod Q_*}
\frac{\log_q|a_{n+1}|}{\log_q\big|a_{n+1}\prod_{i=1}^n a_i^2\,\big|}\;.
\label{eq:middlememb}
\end{align}
Now Theorem \ref{theo:application} follows from Equation
\eqref{eq:premiermemb}, Equation \eqref{eq:secondmemb}, and Equation
\eqref{eq:middlememb}. \cqfd

\bigskip
\noindent
\rem Note that the definition of $\Delta_{Q_*}$ is related to the
choice of one of the ends of $\Ga^0_{Q_*}\bs \TT_{\wh K}$ (the one
corresponding to the geodesic ray in $\TT_{\wh K}$ with vertices
$[\OOOO\times X^{-n}\OOOO]$ for $n\in\NN$).  Other choices of ends
give analogous Diophantine approximation results.

\medskip
{\small \bibliography{../biblio} }

\bigskip
{\small\noindent \begin{tabular}{l} 
Department of Mathematics, University of Illinois Urbana-Champaign\\
1409 W. Green Street, URBANA IL 61801, USA\\
{\it e-mail: jathreya@illinois.edu}
\end{tabular}
\medskip

\noindent \begin{tabular}{l}
D\'epartement de math\'ematique, UMR 8628 CNRS, B\^at.~425\\
Universit\'e Paris-Sud,
91405 ORSAY Cedex, FRANCE\\
{\it e-mail: frederic.paulin@math.u-psud.fr}
\end{tabular}

}

\end{document}